\documentclass[11pt]{article}
\usepackage{amsmath,amsfonts,amssymb,theorem,graphicx,verbatim,fullpage}
\usepackage{hyperref}
\usepackage{pb-diagram}
\newcommand{\CC}{\mathbb C}

\newcommand{\AAA}{\mathbb A}
\newcommand{\PP}{\mathbb P}
\newcommand{\ZZ}{\mathbb Z}

\newcommand{\QQ}{\mathbb Q}

\newcommand{\lam}{\lambda}

\newcommand{\Lam}{\Lambda}
\newcommand{\al}{\alpha}

\newcommand{\gog}{\mathfrak g}

\newcommand{\got}{\mathfrak t}
\newcommand{\gof}{\mathfrak f}

\newcommand{\gop}{\mathfrak p}

\newcommand{\nabp}{\nabla_{\gop}}

\newcommand{\la}{\lambda}
\newcommand{\om}{\omega}
\newcommand{\si}{\sigma}
\newcommand{\SiF}{\Sigma F_4}

\newcommand{\Si}{\Sigma}

\newcommand{\Oh}{\mathcal O}

\newcommand{\Cc}{\mathcal C}

\newcommand{\into}{\hookrightarrow}
\newcommand{\onto}{\twoheadrightarrow}
\newcommand{\ra}{\rightarrow}

\newcommand{\color}[6]{}

\theorembodyfont{\rmfamily}
 \newtheorem{theorem}[subsection]{Theorem}

\newtheorem{example}[subsection]{Example}
\newtheorem{rmk}[subsection]{Remark}
{
\theorembodyfont{\rmfamily}

 \newtheorem{nothing}[subsection]{}
}
 
 \newenvironment{proof}{\paragraph{Proof}}{\par\medskip}

\theorembodyfont{\rmfamily}

\DeclareMathOperator{\fl}{FL}

\DeclareMathOperator{\Gr}{Gr }
\DeclareMathOperator{\lgr}{LGr }
\DeclareMathOperator{\Hom}{Hom}

\DeclareMathOperator{\Proj}{Proj}

\DeclareMathOperator{\Char}{Char}

\newtheorem{thm}{Theorem}[section]
\usepackage[normalem]{ulem}
 \newtheorem{dfn}[thm]{Definition}

\numberwithin{equation}{section}
\makeatletter
\newcommand{\subjclass}[2][2000]{%
  \let\@oldtitle\@title%
  \gdef\@title{\@oldtitle\footnotetext{#1 \emph{Mathematics subject classification.} #2}}%
}
\newcommand{\keywords}[1]{%
  \let\@@oldtitle\@title%
  \gdef\@title{\@@oldtitle\footnotetext{\emph{Key words.} #1.}}%
}
\makeatother
\title{ Polarized 3-folds in a codimension 10 \\ weighted homogeneous \(F_4\) variety }

\date{}
\author{Muhammad Imran Qureshi\\{\sc  Department of Mathematics, SBASSE, LUMS }
\\ {\sc Opp. Sector-U DHA, Lahore, Pakistan.}
\\
{\it Email address}: {\tt i.qureshi@maths.oxon.org}}

\begin{document}

\maketitle

\begin{abstract}
We describe the construction of a codimension 10 weighted homogeneous variety \(w\Si F_4(\mu,u)\) corresponding to the exceptional Lie group \(F_4\) by explicit computation of its graded ring structure. We  give a  formula for the  Hilbert series of the generic weighted \(w\SiF(\mu,u)\)   in terms of  representation theoretic data of \(F_4\). We also  construct some families of polarized 3-folds in codimension 10  whose general member is a weighted complete intersection of  some  \(w\Si F_4(\mu,u)\).\end{abstract}

\section{Introduction}
  We are interested in the study of projective algebraic varieties in terms of graded rings, which are usually Gorenstein in interesting cases.  In lower codimension, one can   describe the  structure  of varieties  by using the standard structure theory of Gorenstein rings: they are defined by a single equation in codimension 1, in codimension 2 they are complete intersections \cite{serre}, and in codimension 3 they are defined by the  \(2m \times 2m  \) pfaffians of a \((2m+1)\times (2m+1)\) skew symmetric  matrix \cite{BE}. The general structure theory for  codimension 4 Gorenstein rings was developed by Reid in \cite{codim4} but, in the words of the author himself, is still some way from any tractable application. Techniques like unprojection have been successfully used in \cite{tom-jerry-1}  to construct  Gorenstein rings in codimension 4 from   rings in lower codimension, but   usually there are    obstructions to using unprojection.  

To construct projective varieties having graded rings in codimension \(>\) 4, weighted homogeneous varieties \(w\Si\) has been used as ambient (key) varieties to construct projective varieties as  weighted complete intersections.  The notion of weighted homogeneous variety  was first introduced by Grojnowski and Corti--Reid in \cite{wg},  a weighted projective analogue of the classical homogeneous variety \(\Si=G/P\); \(G\) a reductive Lie group and \(P\) is a parabolic subgroup of \(G\). They constructed some families of polarized 3-folds which are weighted complete intersections of the weighted orthogonal Grassmannian OGr(5,10), by computing the corresponding Hilbert series and the graded ring structure of \(w\)OGr(5,10). The orthogonal Grassmannian OGr(5,10) is the quotient of the  even orthogonal group O(10,\(\CC)\) by one of its maximal parabolic subgroups, .

 A  formula for the Hilbert series of any weighted homogeneous variety and an algorithmic approach  to compute their defining ideals has been given by Qureshi and Szendr\H oi  in \cite{qs}. This  lead to an efficient approach of calculating graded rings and  required free resolution information for any weighted homogeneous variety, allowing for the construction of projective varieties in higher codimension.      In \cite{qs,qs-ahep}, a weighted $G_2$ variety was used to construct  families of  polarized 3-folds in codimension eight and the weighted Grassmannian $w\Gr(2,6) $  in codimension 6; to construct the polarized varieties in relatively higher codimension.  The weighted Lagrangian Grassmannian \(w\lgr(3,6)\) and   the weighted partial \(A_3\) flag variety \(w\fl_{1,3}\) have been used to construct families of polarized 3-folds in  codimension \(7\) and 9 respectively in \cite{qs2}.

In this article, we construct some families of  polarized 3 dimensional orbifolds, i.e. 3-folds with worst terminal quotient singularities,  in codimension 10. We  explicitly describe  the construction   of a new weighted homogeneous variety \(\left(w\SiF,\Oh_{w\SiF}(1)\right)\):   a homogeneous variety for   the  simple Lie  group $F_4$.   Our $w\SiF$ is a 15 dimensional variety and has an embedding in the weighted projective space $w\PP^{25}$:  a codimension ten embedding. We also give a compact formula for the Hilbert series of a generic weighted homogeneous  \(w\SiF(\mu,u). \)   \\ 
The explicit   graded ring  construction  of \(w\SiF\)  will be given  in terms of generators and relations, by using an algorithmic approach \cite[Appendix A]{qs}.  Then we construct  some  families of polarized 3-folds  as weighted completed intersections of \(w\SiF\) by using  the graded rings and  Hilbert series of \(w\SiF\). The graded rings of these polarized 3-folds  by using the defining equations of \(w\SiF\).
   
    We construct  3-folds \((X,D)\)  polarized by \(\QQ\)-ample Weil divisor    \(D \)  having a  finitely generated graded ring \begin{displaymath}
R(X,D)=\bigoplus_{m \geq 0}H^0(X,mD).
\end{displaymath}
The  surjective morphism from a free graded polynomial ring generated in degree \(w_{i}\) by the variables \(x_i\) \begin{displaymath}
\CC[x_0, \cdots , x_n]\onto R(X,D)
\end{displaymath}
gives an embedding \[i : X = \Proj R(X,D) \into \PP[w_0, \cdots , w_n].\]     The  divisorial sheaf \(\Oh_X(D)\),  a  rank one reflexive sheaf, of \(X\) is  isomorphic to \(\Oh_X(1) = i^*\Oh_{\PP}(1). \)

  We use the natural Pl\"ucker-type embeddings   of  the  corresponding weighted  homogeneous  variety $(w\SiF, \Oh{w\SiF}(1))$ to construct   examples in codimension 10 by taking  quasilinear sections.  We calculate the Hilbert series of a given weighted homogeneous variety to compute the canonical divisor class of \(w\SiF\), corresponding to the  choice of parameters \(\mu\) and \(u\).  Then to  construct a variety with a required canonical class we take  the complete intersection of \(w\SiF\) or of  projective cone(s)  over it, with hypersurfaces  in weighted projective space of appropriate degree. We need the defining ideals  of \(w\SiF\) varieties to understand the type of singularities of their complete intersections. The defining ideals of homogeneous  varieties appear  in~\cite[Sec 1]{rudakov}. We compute the equations by using  the algorithmic approach of DeGraaf \cite{degraaf}.  

 The Section~\ref{notationsec} gives  the required definitions and   conventions used in the rest of the article which includes a recall of the weighted homogeneous varieties, a   formula for  their Hilbert series $P_{w\Sigma}(t)$ and computations of their  defining equations. In Section \ref{f4sec}, we study the structure of weighted homogeneous variety \(w\SiF(\mu,u)\) and its quasilinear sections, leading to  codimension ten polarized varieties.

\subsection*{Acknowledgements}
 I am thankful to  Bal\' azs Szendr\H oi     for the  helpful discussions and an anonymous referee for the suggestions to  improve an earlier version of this article.  This research is supported by the Lahore university of Management Sciences's (LUMS)  faculty startup research grant.
\section{Definitions and  notations }

\label{notationsec}

We work over a field complex numbers \(\CC\). A    polarized variety is a pair  \((X,D);\)   ~\(X\) is a normal projective  algebraic variety and~\(D\) is  
a \(\QQ\)-ample  Weil divisor on $X$. All our varieties  are well-formed and quasi-smooth, embedded in some weighted projective space. We use the  standard notations
$\PP[w_0,w_1,\cdots,w_n]$ or \(\PP^n[w_i]\)  to denote the   weighted projective space.

A polarized variety \(X\subset \PP^n[w_i]\) of codimension \(l\)  is called well-formed, 
if no \(n\) of \(w_0,\cdots,w_n\) have a common factor  and singularities of \(X\)only  appear in  codimension greater than \( l+1\). If the affine cone \(\widetilde X \subset \AAA^{n+1}\) over 
\(X\) is smooth outside the origin  then \(X\) is called quasi-smooth. If \(X\) is quasi-smooth, then it contain no other singularities than those coming from the weights of the embedding  \(\PP^n[w_i]\), called quotient singularities. We assume that  restriction of the tautological ample divisor 
$\Oh_{\PP}(1)$  provides the polarization. 

The Hilbert series of  \((X,D)\) is given by  
\[ P_{(X,D)}(t)=\sum_{m\geq 0}\dim H^0\left(X,\Oh_X(mD)\right) t^{m}.
\]
 Appropriate vanishing theorems and Riemann--Roch 
formulas can be used together  to compute $h^0(X,mD)=\dim H^0(X, mD)$ 
in most of the  cases. We will  write $P_X(t)$ for the Hilbert series if no confusion can arise. A polarized   3-fold \((X,D)\) is a three dimensional Gorenstein, normal, projective  algebraic variety with at worst terminal  cyclic quotient 
singularities, consisting of isolated orbifold points of \(X\).  A projective variety \(X\) is called Gorenstein if \begin{enumerate}
\item \(H^i(X,\Oh_X(mD))=0 \) for all \(m\ge 0\) and \(0<i<\dim(X)\).
\item \(K_X\sim kD\) for some some integer \(k\).

\end{enumerate} A cyclic quotient singularity \(Q\) of type \(\dfrac1r (a_1,\ldots,a_n)\) is the quotient \(\AAA^n /\ZZ_r\) given by \[\zeta:(x_1,\ldots,x_n)\mapsto \left(\zeta^{a_1}x_1,\ldots,\zeta^{a_n}x_n\right ), \]  where \(\zeta\) is a primitive \(r\)-th root of unity. It is called isolated if \( \gcd(r,a_i)=1\) for all \(1\le i\le n\) and terminal if  \[\frac1r\sum_{j=1}^{n} \overline{d a_j} >1 \text{ for }k=1,\ldots,r-1,\] where \(\overline{da_j}\) is the smallest residue modulo \(r\), see \cite{YPG}.  

Now we recall some background from representation theory which is primarily from  \cite{harris}. Let  \(G \)  be a linear algebraic group, \(P\) be a parabolic subgroup of \(G\) and \(T\) be the maximal torus inside \(G\) then  \(T\subset  P \subset G\). Then \(\got\subset \gop \subset \gog\) are the inclusions of the  corresponding Lie algebras. The quotient \(\Si=G/P\) is 
called a homogeneous variety or generalized flag variety.
 Let $\Lam_W={\rm Hom}(T, \CC^*)$ denote the  lattice of weights.
Let \(V\) be an irreducible \(G\)-representation then \(V\) has a decomposition 
\[V=\bigoplus_{\al\in\Lambda_W} V_{\al}\]
into eigen\-spaces of  \(T\)-action;
the \(\al\)'s with non-trivial $V_\alpha$ are called the weights of the representation \(V\). If \(V=\gog\) then \(V\) is called the adjoint representation and  the non-zero weights are called the roots \(\nabla \) of \(\gog\). Each element of \(\nabla\) can be written as  a strictly positive or negative linear combination of the subset \(\nabla_0\), called the simple roots of \(\gog.\) The set of roots  \( \nabla\) has a decomposition into the set of positive and negative roots  \(\nabla=\nabla_+\sqcup \nabla_ - \).

Let   $V_\chi$ be  the irreducible representation of \(G\) with highest (dominant)
weight   \(\chi\) with respect to certain partial order on the set of weights and  \(\nabla(V)\) represents the set of weights of  \(V_{ \chi}\). Then the  set of parabolic subgroups (up to conjugacy) of \(G\) are in 1--1 correspondence with the irreducible highest weight representations of \(G\).  The character  \(\chi: T\to \CC^* \) give rise to a very  ample line bundle \(\mathcal{L}_\chi\) on the homogeneous variety \(\Sigma\). Then by Borel--Bott--Weil theorem \(H^0(\Si,\mathcal L_\chi)=V_\chi\) and we get the embedding   \[\Si=G/P_\chi\into\PP V_\chi,\] where     
$P = P_\chi$ is the parabolic corresponding to the subset of simple roots which are orthogonal to   \(\chi\) in the character lattice \(\Lam_W\).
The dimension of \(V_{\chi}\) can be calculated by using the Weyl's dimension formula: \begin{equation}\label{wdf}
\dim (V_\chi)=\prod_{\al \in \nabla_+ }\dfrac{(\chi+\rho,\al)}{(\rho,\al)},
\end{equation} 
where the Weyl vector  \(\rho\) is half the sum of the positive roots  and \((,)\) is the Killing form on the Lie algebra \(\gog\). The full character of \(V_{\chi}\) can be computed by using the Weyl character formula:
\begin{equation}\label{WCF}\Char(V_\chi)=\dfrac{\displaystyle\sum_{\si \in W} (-1)^\si t^{\si(\chi+\rho)}}{\sum_{\si \in W} (-1)^\si t^{\si(\rho)}}=\sum_{\chi_i\in \nabla(V)}\dim(V_{\chi{_i}})t^{\chi_i},\end{equation} where \(W\) is the Weyl group of the root system of \(\gog \) and \(\nabla(V)\) is the set of weights of the representation \(V\). If $\si$ is a product of   an even  number of simple reflections in  $W$ then $(-1)^\si=1  $ and -1 otherwise.    
Let $\Lam_W^* = \Hom(\CC^*,T)$ denote  the  lattice of 1-parameter subgroups of \(G\), then we have get a perfect pairing \(<,>:\Lam_W\times \Lam_W^*\ra\ZZ. \)  We take an element  \(\mu\) in \(\Lam_W^* \) and a positive integer \(u \) such 
that
 \begin{equation}
\left<\si\chi,\mu\right>+u >0,
\label {weights} 
\end{equation}
for all  elements \(\si\) of the Weyl group \(W\). The   inequality \eqref{weights} makes sure that  all the weights on the weighted projective space containing \(w\Si(\mu,u)\) are positive.

\dfn \cite{wg}  
Let \(\Si\) be a homogeneous variety.  
Take the affine cone $ \widetilde{\Sigma} \subset  V_{\chi}$ of the 
embedding $\Sigma \hookrightarrow \PP V_{\chi}$
then the invariant part of  the following $ \CC^*$-action  on $V_{\chi}\backslash\{0\}$ 
\[ 
(\varepsilon \in \CC^*) \mapsto ( v \mapsto \varepsilon^u(\mu(\varepsilon)\circ v)),
\] is called the weighted homogeneous  variety embedded  in \(\PP[<\chi_i,\mu>+u]\), where \(\chi_is \) are the weights of the presentation \(V_\chi\).  We denote 
 this variety by \(w\Si(\mu,u)\) as the weights depends on the paramenters \(\mu\) and \(u\) or simply $w\Sigma$, if no confusion can arise. 
 
The  following formula from \cite[Thm. 3.1]{qs} can be used to figure out the    Hilbert series  of any  weighted homogeneous  variety.
\begin{equation}
P_{w\Si}(t)=\dfrac{\displaystyle\sum_{\si\in W}(-1)^\si \dfrac{t^{\left<\si\rho, \mu\right>}}{(1- t^{\left<\si\chi,\mu\right>+u})}}{\displaystyle\sum_{\si\in W}(-1)^\si t^{\left<\si\rho, \mu\right>}},
\label{whhs}
\end{equation}

\begin{rmk}  The Hilbert series \eqref{whhs} simplifies to the   expression
\begin{equation}
P_{w\Si}(t)=\dfrac{\ N(t)}{\displaystyle \prod_{\chi_i \in \nabla{(V_\chi)}}(1-t^{<\chi_i,\mu>+u})},
\label{reducedhs}
\end{equation} by using the standard 
Hilbert--Serre theorem \cite[Theorem 11.1]{atiyah}, where \(\chi_i\) are the weights of the representation \(V_\chi\).

\end{rmk}

The defining ideal \(I =\left<Q\right>\)  of a homogeneous variety \(\Si=G/P \into \PP V_{\chi}\)  is always generated by  quadrics \cite[2.1]{rudakov}. The  second symmetric power of  \(V^*_{\chi}\) has a  decomposition   
\begin{displaymath}
S^2 (V_\chi^*)=V_{2\nu}\oplus V_1 \oplus\cdots \oplus V_n
\end{displaymath}  
into irreducible   \(G\)-representations, where  \(\nu\) is  the highest weight of  \(V^*_{\chi}\). The generators of the  subspace  \(Q \subset S^2 V^*_{\chi}  \)
 consisting of all the summands except~\(V_{2\nu}\), gives  the defining equations of \(\Si\).    

% \subsection*{Constructing polarized varieties}
% We start by constructing  some key variety \(w\Si \) (weighted flag variety) embedded into some weighted projective  space.  Then we find the Hilbert series of the given variety, which gives us some information about the graded ring \(R(w\Si,D)\).  Under suitable  conditions, we can compute   the canonical divisor class of  \(w \Si\). Then we take  quasilinear sections (general hypersurfaces of appropriate degree in weighted projective space) of \(w\Si\) or of  projective cones  over it, to get a  variety with the desired canonical or anticanonical class. Then we study  different aspects, such as  singularities, well-formedness and quasi-smoothness of the resulting variety to establish the existence of an appropriate model of the variety. At the end, 
% using the orbifold 
% Riemann--Roch formula of~\cite[Section 3]{anita}, we  compute the invariants of our polarized variety \((X,D)\)    from the first few values of $h^0(nD)$, and verify 
% that the same Hilbert series can be recovered.
%  More  details can be found in \cite{qs}. 
\section{Weighted  homogeneous $F_4$ varieties}
\label{f4sec}
In this section, we recall the representation theory of the Lie group \(F_4\) and how we can use it to construct the weighted homogeneous variety \(w\SiF \). We also compute  a  formula for the Hilbert series of \(w\SiF \). At the end, we construct some families of 3-folds of general type as weighted complete intersection of some   \(w\SiF.\)
\subsection{Generalities}\label{basicsf4}In this section, we recall the  basic algebraic structure of the Lie group \(F_{4}\)  and its corresponding Lie Algebra \(\gof_4\) which can mostly be found in \cite{bourbaki46}. Let \(G\) be the simple and simply connected   Lie group of type \(F_4\)  with Lie algebra  \(\gog=\gof_4\). This  is one of the five exceptional  Lie groups.  The rank of the Lie algebra \(\gof_4 \) is 4, which is the dimension of the Cartan subalgebra \(\got\) of \(\gof_4.\) The weight lattice of  \(\gof_4\) is  a rank four lattice\[\Lam_W=\left<e_1,e_2,e_3,e_4\right> .\]  The  fundamental weights of  $G$, which form a set of generators of the dominant Weyl chamber,  are given by \[\om_1=e_1+e_{2},\text{\;\;\;}\om_2=e_1+e_2+e_3,\text{\;\;\;}\om_3=\dfrac{1}{2}(3e_1+e_2+e_3+e_4),\text{\;\; and }\;\; \om_4=e_{1}.\]  
 The  simple roots  of the  root system of \(\gof_4\) are \[\al_1=2\om_{1}-\om_{2}=e_2-e_3,\text{\;\;\;\;\;\;\;\;\;\;\;\;\;\;\;\;\;\;\;\;\;\;\;\;\;\;\;\;\;\;\;\;\;\;\;\;}\al_2=-\om_{1}-2\om_3+2\om_{2}=e_3-e_4\]\[\al_3=-\om_4+2\om_3-\om_2=e_4,\text{\;\; and \;\;\;\;\;\;\;\;\;\;\;\;\;\;}\al_4=2\om_{4}-\om_3=\dfrac{1}{2}(e_1-e_2-e_3-e_4).\]  The Weyl group \( W\)  of the root system  \(F_4\) is generated by the   simple reflections \(s_{\al_i},\; 1 \leq i \leq 4\), and  is the   symmetry group of the convex regular 4-polytope, known as the 24-cell in convex geometry. The order of \(W\) is 1152 and has  a finite presentation  given as follows. \[W=\left<s_{\al_i}\vert s_{\al_i}^2=(s_{\al_1}s_{\al_2})^3 = (s_{\al_1} s_{\al_3})^2 = (s_{\al_1}s_{\al_4})^2 = (s_{\al_2}s_{\al_3})^4 = (s_{\al_2}s_{\al_4})^2 = (s_{\al_3}s_{\al_4})^3 =1\right>\]The Weyl vector \(\rho\) is; 
\begin{displaymath}\rho=8\al_1+15\al_2+21\al_3+11\al_4=\dfrac{1}{2}(11e_1+5e_2+3e_3+e_4).\end{displaymath}
Consider   the  irreducible fundamental  representation  \(V_\chi\) of the Lie group \(F_4 \) with highest weight \(\chi=\om_4=e_1\). Then  the dimension of the irreducible representation \(V_{\chi}\) is 26 by using  the  Weyl dimension formula \eqref{wdf}. Twenty-four  of the weights appear with multiplicity one and the zero weight space appears with multiplicity two. This can be easily  figured out by using the implementation of Weyl character formula \eqref{WCF} in some computer algebra system like  SAGE \cite{sage} or LiE \cite{LiE}. Sixteen of the  24 non-zero weights appear with multiplicity one are given by 
\begin{equation}\label{wt1}\left\{\dfrac12\left((-1)^ie_1+(-1)^je_2+(-1)^ke_3+(-1)^le_4\right):0\le i,j,k,l\le 1\right\}\end{equation} and eight of them are  \begin{equation}\label{wt2}\left\{(-1)^ie_j:0\le i\le 1,1\le j\le 4\right\}.\end{equation}

% \renewcommand{\arraystretch}{1.5}\label{weightsf4}\begin{eqnarray}\left.
% \begin{array}{lll}  \chi=\chi_1= e_1 & \chi_2=\frac{1}{2} (-e_1-e_2-e_3-e_4), & \chi_3=-e_1\\
%                     \chi_4=\frac{1}{2} (e_1+e_2+e_3+e_4)& \chi_5=\frac{1}{2} (e_1-e_2-e_3-e_4)& \chi_6=\frac{1}{2} (e_1+e_2+e_3-e_4)\\
%                     \chi_7=\frac{1}{2} (e_1-e_2-e_3+e_4)&\chi_8=\frac{1}{2} (e_1+e_2-e_3+e_4),&\chi_9=\frac{1}{2} (e_1-e_2+e_3-e_4)\\
%                     \chi_{10}=\frac{1}{2} (e_1-e_2+e_3+e_4)&\chi_{11}=\frac{1}{2} (e_1+e_2-e_3-e_4)&\chi_{12}=\frac{1}{2} (-e_1+e_2+e_3+e_4)\\
%                     \chi_{13}=e_4&\chi_{14}=e_3&\chi_{15}=e_2\\
%                     \chi_{16}=\frac{1}{2} (-e_1+e_2+e_3-e_4)&\chi_{17}=-e_4&\chi_{18}=\frac{1}{2} (-e_1+e_2-e_3+e_4)\\
%                     \chi_{19}=-e_3&\chi_{20}=\frac{1}{2} (-e_1-e_2+e_3+e_4)&\chi_{21}=-e_2\\
%                     \chi_{22}=\frac{1}{2} (-e_1+e_2-e_3-e_4)&\chi_{23}=\frac{1}{2} (-e_1-e_2+e_3-e_4)&\chi_{24}=\frac{1}{2} (-e_1-e_2-e_3+e_4)\\
%                     &\chi_{25}=\chi_{26}=0  \\   
%  \end{array}\right\}\end{eqnarray}

 The Lie algebra $ \gof_4 $ is the  adjoint representation of the Lie group \(F_4\). In fact it is the irreducible representation with highest weight \(\om_1=e_1+e_2\). One can use       the  formula \eqref{wdf} to show that it is 52 dimensional.  The  simple roots $ \al_2,\al_3 \text{ and }\al_4 $ are  orthogonal to the dominant weight $\chi=\om_4$ with respect to the  the Killing form.  Therefore, the  parabolic subalgebra is   \[\gop_{\chi}= \bigoplus\left(\got \bigoplus_{\al \in \nabla_+}\gog_{\al}\bigoplus_{\al \in {\nabp}} \gog_{-\al}\right),\]where $\nabp$  is the subset of $ \nabla_+$  consisting of the terms involving only  $\al_2,\al_3 \text{ and }\al_4$. The  parabolic subalgebra \(\gop_{\chi} \) is 37 dimensional. Consider the  quotient of the Lie group \(G\) by the parabolic subgroup \(P_{\chi},\) corresponding to the subalgebra \(\gop_{\chi},\) then   \(\SiF=G/P_{\chi}\) is a homogeneous variety embedded in the projectivization  of the \(G\)-representation \(V_{\chi}\). The dimension of the homogeneous variety \(\SiF \text{ is } 52-37=15\)  and therefore we get a codimension 10 embedding  \[\SiF \into \PP^{25}[V_{\chi}].\] We denote the given homogeneous variety by \(w\SiF\). Let  \[\Lam_W^*=\left<f_1,f_2,f_3,f_4\right>\] The weighted version of the homogenous variety \(\Si\) is obtained by choosing \[\mu= \displaystyle\sum_{i=1}^4 a_i f_i  \in \Lam_W^* \] and a positive integer  \(u\) such that \(\left<s_{\al}\cdot \chi,\mu\right>+u >0\) for all the elements \(s_\al\) of   \(W\).  Then  for each such choice of \(\mu\) and \(u\) we  have the embedding of  weighted homogeneous variety
\begin{equation}\label{weights}w\Si^{15}(\mu,u)\into w\PP^{25} V_{\chi}[\left<\chi_i,\mu\right>+u ],\end{equation} where  \(\chi_i\)'s are the weights of the representation \(V_{\chi}\), \(i=1,\cdots, 26\). The element \(\mu \) of the dual lattice is usually represented as a vector \(\mu=(\uline {a_i}).\)

\subsection{Hilbert series of weighted  $w\SiF(\mu,u)$ } 
We   use  computer algebra systems SAGE and Mathematica together to compute the  following Hilbert series formula of \(w\SiF\). The SAGE has built in representation theoretic data which has been used to compute the numerator and denominator of the formula \eqref{whhs}. A  SAGE code  is given in the  Appendix \ref{sagecode} to compute the Hilbert series of \(w\SiF\) which in principal can be modified to compute the Hilbert series of any weighted homogeneous variety.  The further simplification has been performed by using Mathematica.  
\thm Consider the symmetric group \(S_2\) with  \(\si(a_i)=a_i\) if \(\si\) is an  even and $\si(a_i)=-a_i$ if \(\si\) is an  odd permutation, where \(a_i \in \mu\). Then  the Hilbert series of the \(w\SiF \)  has  the following  compact form.
\begin{equation}\label{hsf4} 
P_{w\SiF}(t)=\displaystyle\dfrac{1-\displaystyle\sum_{k=2}^3(-1)^k P_{k-1}(t)\left( t^{ku}+t^{(15-k)u}\right)+\displaystyle\sum_{k=5}^7(-1)^k P_{k-2}(t)\left( t^{ku}+t^{(15-k)u}\right) -t^{15u}}{\displaystyle\prod_{\chi_i \in\nabla(V_\chi)} \left( 1-t^{<\chi_i,\mu>+u} \right)},\end{equation}

where 
 \[\begin{array}{cc}
P_1(t)&=\displaystyle\sum_{\sigma\in S_2}\displaystyle\sum_{i=1}^4 t^{\si(a_i)}+\displaystyle\sum_{0 \leq(i,j,k,l)\leq 1}t^{\frac12\left((-1)^ia_1+(-1)^ja_2+(-1)^ka_3+(-1)^la_4\right)}+3\end{array},
\]
 \[ \begin{array}{cc}
 P_2(t)&=\displaystyle\sum_{\sigma\in S_2}\displaystyle\sum_{i=1}^4 2t^{\si(a_i)}+\displaystyle\sum_{0 \leq(i,j,k,l)\leq 1}2t^{\frac12\left((-1)^ia_1+(-1)^ja_2+(-1)^ka_3+(-1)^la_4\right)}\\
&
+  \displaystyle\sum_{0\le (i,j)\le 1}\left(\sum_{1\le m<n\le 4} t^{(-1)^i a_m+(-1)^ja_n}  \right)+6
\end{array}, \] 
\[ \begin{array}{cc}
 P_3(t)&=\displaystyle\sum_{0 \leq(i,j,k,l)\leq 1}\left(7t^{\frac12\left((-1)^ia_1+(-1)^ja_2+(-1)^ka_3+(-1)^la_4\right)}+\sum_{m=1}^4t^{a_m} t^{\frac12\left((-1)^ia_1+(-1)^ja_2+(-1)^ka_3+(-1)^la_4\right)}\right)\\
&
+  \displaystyle\sum_{0\le (i,j)\le 1}\left(\sum_{1\le m<n\le 4} 3t^{(-1)^i a_m+(-1)^ja_n}  \right)+\sum_{0\le (i,j,k)\le 1}\left(\sum_{1\le l< m<n\le 4} t^{(-1)^i a_l+(-1)^ja_m+(-1)^ka_n}  \right)\\&+\displaystyle\sum_{\sigma\in S_2}\displaystyle\sum_{i=1}^4 7t^{\si(a_i)}+15
\end{array} ,\]
\[ \begin{array}{cc}
 P_4(t)&=\displaystyle\sum_{0 \leq(i,j,k,l)\leq 1}\left(12t^{\frac12\left((-1)^ia_1+(-1)^ja_2+(-1)^ka_3+(-1)^la_4\right)}+\sum_{m=1}^42t^{a_m} t^{\frac12\left((-1)^ia_1+(-1)^ja_2+(-1)^ka_3+(-1)^la_4\right)}\right)\\
&
+  \displaystyle\sum_{0\le (i,j)\le 1}\left(\sum_{1\le m<n\le 4} 5t^{(-1)^i a_m+(-1)^ja_n}  \right)+\sum_{0\le (i,j,k)\le 1}\left(\sum_{1\le l< m<n\le 4} 2t^{(-1)^i a_l+(-1)^ja_m+(-1)^ka_n}  \right)\\&+\displaystyle\sum_{\sigma\in S_2}\displaystyle\sum_{i=1}^4\left( 12t^{\si(a_i)}+t^{\si(2a_i)}\right)+\displaystyle\sum_{0 \leq(i,j,k,l)\leq 1}t^{(-1)^ia_1+(-1)^ja_2+(-1)^ka_3+(-1)^la_4}+26
\end{array} ,\]and 
\[ \begin{array}{cc}
 P_5(t)&=\displaystyle\sum_{0 \leq(i,j,k,l)\leq 1}\left(6t^{\frac12\left((-1)^ia_1+(-1)^ja_2+(-1)^ka_3+(-1)^la_4\right)}+\sum_{m=1}^4t^{a_m} t^{\frac12\left((-1)^ia_1+(-1)^ja_2+(-1)^ka_3+(-1)^la_4\right)}\right)\\
&
+  \displaystyle\sum_{0\le (i,j)\le 1}\left( \sum_{1\le m<n\le 4}3t^{(-1)^i a_m+(-1)^ja_n}  \right)+\sum_{0\le (i,j,k)\le 1}\left(\sum_{1\le l< m<n\le 4} t^{(-1)^i a_l+(-1)^ja_m+(-1)^ka_n}  \right)\\&+\displaystyle\sum_{\sigma\in S_2}\displaystyle\sum_{i=1}^4\left( 6t^{\si(a_i)}+t^{\si(2a_i)}\right)+\displaystyle\sum_{0 \leq(i,j,k,l)\leq 1}t^{(-1)^ia_1+(-1)^ja_2+(-1)^ka_3+(-1)^la_4}+15
\end{array} .\] Moreover, if 
\(w\Si\) is  well-formed  then the  canonical line bundle \(K_{w\Si}=\Oh_{w\Si}(-11u).\)
 
\proof  We first calculate the  orbit of the  weight \(\chi\) under the action of the Weyl group \(W\).   From the representation theory since the image of the non-zero weight is a non-zero weight under he action of Weyl group,  only non-zero weight    appear in the orbit of \(\chi\). By using SAGE we compute that all of the  24 non-zero weight appear in the orbit \(W\chi\).  We evaluate the formula  \eqref{whhs} for  \(W,\rho, \mu\) and \(\chi,\) as given in  Section \ref{basicsf4} to compute \(P_{w\SiF}(t); \) by using SAGE code in Appendix \ref{sagecode}. Then we perform simplification in Mathematica to obtain the  below form of the Hilbert series of \(w\SiF\).
\begin{equation}\label{hs:f4}P_{w\SiF}(t)=\dfrac{1+2t-\displaystyle \sum_{i=1}^{24}t^{\left<\chi_i,\mu\right>+2u}+\cdots +2t^{12u}+t^{13u}}{\displaystyle\prod_{i=1}^{24}\left(1-t^{\left<\chi_i,\mu\right>+u}\right)},\end{equation} where \(\chi_i\) are the  collection of weights of \(V_{\chi }\) given by \eqref{wt1} and \eqref{wt2}. The zero weight spaces do not appear in the orbit under the action of the Weyl group of \(F_4\). Therefore the full expression of type \eqref{reducedhs} for the Hilbert Series of \(w\SiF\) is obtained  obtained by multiplying and dividing the equation \eqref{hs:f4}  by \((1-t^u)^2\), which represents the  zero weight spaces in the representation \(V_{\chi}\) as \(\chi=(\underline 0) \) in \eqref{weights}. This gives us the compact form \eqref{hsf4} of the Hilbert series. 
The sum of the weights on \(\PP[\left<\chi_i,\mu\right>+u]\) is 26\(u\) and the adjunction number  is \(15u\). Therefore, if  \(w\SiF\) is well-formed (normal), then the canonical divisor class  is \(K_{w\SiF}=\Oh(15u-26u)=\Oh(-11u).\) \subsection{  Families of 3-folds in codimension 10 }   In this section, we construct some  families of 3-folds  of general type with at worst terminal quotient singularities  with  the  canonical divisor  class \(K_X=kD, k \geq 5\).  We also tried to construct  some families of  Calabi--Yau and \(\QQ\)-Fano threefolds   in \(w\SiF\) by using the algorithmic approach of \cite{QJSC} but the computer search  was not successful.     We   show that  a threefold  linear section  of the straight homogeneous variety \(\Si\)  is a smooth  canonical 3-fold. The search for canonical 3-folds with non-trivial weights was also unsuccessful.
\begin{example} \label{f4st}  We consider the Hilbert series  straight \(F_{4}\)  variety \(\SiF\), which corresponds to choice of parameters \(\mu=(0,0,0,0)\) and \(u=1.\) The Hilbert series of \(\SiF\) is given by
\[P_{\SiF}(t)=\dfrac{1-27t^2+78t^3-351t^5+ \cdots-351t^{10}+78t^{12}-27t^{13}+t^{15}}{(1-t)^{26}}.\] The canonical divisor is \(K_{\SiF}= \Oh(-11),\) since \(\SiF\) is a  smooth projective variety. Consider the complete intersection  \[X=\SiF \cap \left( \displaystyle\cap_{i=1}^{12} H_i\right),\] where \(H_i\) represents  a  general hyperplane of \(\PP^{25}\). Then \(X\) is a smooth 3-fold with \[K_{X}=\Oh(-11+12)=\Oh(1).\] By using the Hilbert series we  can compute the degree of \(X\) to be 78. Thus \((X,K_{X})\) is a smooth  canonical threefold of general type with  \((K_X)^3=78\).\end{example} 
\rmk In fact we have a ladder of polarized varieties coming from  Example \ref{f4st}. We can easily see that a 5-fold section   \(V\) of \(\SiF\) is a Gorenstein Fano 5-fold polarised by \(-K_V\).  A general member \(Y\in \left|-K_V\right|\) is a Calabi--Yau 4-fold polarized by \(D=-K_V|_{Y}\), and  \(X\) is a general member of  \(\left|D\right |\) polarized by \(K_X=D|_{X}\). Thus we have a ladder of  varieties  \[X\subset Y\subset V\subset \Si.\]
\begin{example}\label{nwf:f4}
We consider the  case  where  ambient weighted homogeneous variety is not well-formed. 
 \begin{itemize}
\item Input parameters: $ \mu=(0,0,0,0) $, $ u=2 $
\item Embedding: $ w\SiF \subset \PP^{25}[\uline 2]  $, with each variable \(x_1,\ldots,x_{26}\) has  weight 2. 

\item Hilbert numerator: $1-27t^4-78 t^6-351 t^{10}+\cdots -351 t^{20}+78 t^{24}-27t^{28}+t^{30}$

\end{itemize} 
To construct a 3-fold of general type we take a triple  projective cone  over \(w\Si\) which makes the given variety well-formed and we have an embedding of an 18-dimensional variety
\[\Cc^3 w\SiF\into \PP^{28}[1^3,2^{26}].\] 
with \(K_{\Cc^3 w\SiF}=\Oh(-25)\); taking each projective cone adds 1 to the dimension and -1 to the canonical weight or index of   $w\SiF$. Then we take the following weighted complete intersection \(\Cc^3w\SiF\) with 15 generic forms \(Q_i\) of degree 2 to get  \[\displaystyle X=\Cc^3w\SiF \cap \left(\cap _{i=1}^{15}Q_{i}\right)\into \PP^{13}[1^3,2^{11}]\] with \(K_X=\Oh_X(-25+2(15))=\Oh_X(5)\) by the adjunction formula. The three new variables  $y_1,y_2, y_3,$ of degree 1 appear in the 15  degree two quadratics form. Each quadric \(Q_i\) replaces a variable \(x_{i} \) of weight 2 with a  from of degree two in the equations of \(X.\) Without any loss of  generality we assume that  \begin{equation}\label{quadrics}x_i:=Q_i(y_{j}s,x_ks)  \text{ for } 12\le i\le 26, 1\le j\le 3,1\le k\le 11;\end{equation} so  we get the embedding     \[X\into\PP[y_1,y_2,y_3,x_1,\ldots,x_{11}].\] The base locus of the linear system $|\Oh(2)|$ of quadrics is empty, due to the weights of  the embedding of \(X\). Thus  the only   singularities of  \(X\) may  occur due to the  non-trivial weights of the embedding.

The locus of weights 2 variables in \(w\SiF \) basically defines the whole \(15 \)-dimensional ambient variety. Then the complete intersection with  15 generic  hypersurfaces of degree 2,  linear equations  in weight 2 variables, given by    \[X_0:=X\cap\{{y_1=y_2=y_3=0\}\into\PP^{10}},\]   is  0-dimensional. Now from the Example \ref{f4st} we conclude that the degree of  \( X_{0}\) is 78. Since \(X_0\) is irreducible and reduced, so it consists of     78 distinct points. This can also be established by using the computer algebra on a specific  example.  Now to determine the local type of each singular point,  we show that on each affine piece of the \(\PP^{10}[x_i]\), the weight 1 variables \(y_1, y_2\) and $y_3$  are  three local variables, so  each  point is locally of  type \(\dfrac12(1,1,1)\).  Now for  affine piece \(x_1\ne 0  \),  we have
 
 \renewcommand*{\arraystretch}{1.5}
$\begin{array}{cccc}
Q_{25}=x_{4}+\cdots, &Q_{12}=x_5+\cdots,&Q_{16}=x_6+\cdots,&Q_{18}=x_7+\cdots\\
Q_{22}=x_{8}+\cdots, &Q_{23}=x_{9}+\cdots,&Q_{24}=x_{10}+\cdots,&Q_{20}=x_{11}+\cdots
\end{array}.$

 By  using the implicit function theorem; from  the equations A.1, A.2, A.3, A.4, A.5, A.6, A.7, A.8, A.9 and A.10 from Appendix \ref{equations},  we can remove the variables \(x_{4}, x_{5},x_{6},x_{7},x_{8},x_{9},x_{10},x_{11},x_{2}\) and $x_{3}$ respectively. Thus \(y_1,y_2, y_3\) are local variables near \(x_1\ne 0\). So the cyclic group \(\ZZ_2 \) acts by \[\zeta:(y_1,y_2,y_3)\mapsto (\zeta y_1,\zeta y_2,\zeta y_3)\] to give the quotient  singularity of type \(\dfrac12(1,1,1)\). We can perform the similar calculation on each affine piece of  \(X\) to show that each point is a singular point of type \(\dfrac12(1,1,1)\).  Thus \(X\) is a family of  well-formed and quasi-smooth polarized 3-folds in the codimension 10 \(w\SiF\). Moreover, the degree of the polarizing divisor \(D^3\) can easily be computed to be 39 by using the Hilbert series of \(X\).  Thus the  degree of the canonical divisor \[(K_X)^3=(5\cdot D)^3=125\cdot39=4875. \]  
\end{example}
\rmk  In fact, we can  construct a finite number of families  of smooth 3-folds of general type  by taking further projective cones over \(w\Si(\uline0,2)\) and taking   3-fold linear sections \(X\) of relatively higher indices \(k\) of \(K_X=\Oh_X(k)\).

\begin{theorem}{ Let  \(w\Si(\uline0,2)\) be the weighted \(F_4\)-variety corresponding to the  choice of parameters \(\mu=(0,0,0,0)\text{ and }u=2 \). Then for each \(k=6,\cdots, 16\), we have a family of smooth  polarized 3-fold $(X_k,D_k)$ of general type of index \(k\), i.e.  \(K_{X_k}=\Oh_{X_k}(kD_k)\) which is the intersection of \(\Cc^{k-2}w\Si(\uline0,2) \) with \((k+10)\) general quadric hypersurfaces of \(\PP^{k+23}\), with the following properties.}
\begin{itemize} {
\item Degree of the embedding: \(D_k^3=39\cdot(2)^{k-5}\)
\item Degree of canonical divisor: \((K_{X_k})^3=(k\cdot D_k)^3\)
\item Weights of the embedding: \(X_k\into \PP^{13}[1^{k-2},2^{16-k}]\)}
\end{itemize} 
\end{theorem}
\proof We start with  the construction of a 3-fold with  \(k=6\). Following the Example \ref{nwf:f4}, if we take a  further projective cone over \( w\Si(\uline0,2) \) then  we add \(-1\) to the  canonical class of  \(w\Si(\uline0,2)\) and  increase its dimension by one as well. Thus by taking four cones over \(w\Si(\uline0,2)  \) and taking a weighted complete intersection with 16 generic forms of degree 2, we get a  3-fold \[X=\Cc^4w\Si\cap(\cap_{i=1}^{16}Q_i)\into \PP^{13}[1^4,2^{10}]\]such that      \(K_X=\Oh(6)\). An extra quadric section will also change the degree \(D^3\) of the new threefold to be \(D^3=39\cdot2=78. \) The degree of the canonical class will be 
\[(K_{X})^3=(6\cdot D)^3=216\cdot78=16848.\] Since we are taking one more quadric section, the singularities of \(X\); thus \[X\cap \{\textrm{ locus of degree 2 variables}\}=\emptyset.\]       
So we get a smooth 3-fold of general type of index 6 in codimension 10. The rest of the cases follow, by recursively using the same line of reasoning for \( k=7,\cdots, 16.\)  
\rmk Recently, an algorithmic approach has been developed to find lists of isolated orbifolds in a given weighted homogeneous variety, see \cite{formats,QJSC}. A more complete list of varieties in codimension 10 \(w\SiF\) varieties will appear elsewhere \cite{BKQ}. 
\appendix \section{Equations of $F_4$ homogeneous variety}\label{equations} We  compute the    equations of the homogeneous variety \(\Si F_4\) by the GAP4 code given in Appendix of  \cite{qs}. We compute the decomposition of the 2nd symmetric power   \(S^2( V_\la^*)\) of the  representation dual to \(V_\la\),  into its direct summands as free modules over  \(\gof_4\).  Now  \(S^2V_\la^*\) is a 351-dimensional vector space and we obtain the following decomposition into the subrepresentations of \(\gof_4\)  \[Z=S^2V_{\lam}^*=V_{1}\oplus V_{e_1} \oplus  V_{e_2},\]
having  dimensions  324, 26, and 1  respectively. The set of  defining equations of \(\Si F_{4}\) is the union of   the basis of the linear subspaces of \(V_{e_1}\) and \(V_{e_2}\).  $$I=\left<Q\right> =\left<V_{e_1} \cup  V_{e_2}\right> \subset S^2V_{\lam}^*.$$  Therefore the homogeneous ideal of \(\Si\) is defined by 27 quadratic equations given below.  
 
\begin{eqnarray} \label{eqf4}
&&\textstyle  x_{1}x_{25}-\frac{1}{2}x_{1}x_{26}-\frac{3}{2}x_{4}x_{5}-\frac{3}{2}x_{6}x_{7}-\frac{3}{2}x_{8}x_{9}+\frac{3}{2}x_{10}x_{11} \\ 
&&\textstyle  x_{1}x_{12}-\frac{1}{3}x_{4}x_{25}-\frac{1}{3}x_{4}x_{26}-x_{6}x_{13}-x_{8}x_{14}+x_{10}x_{15} \\
 &&\textstyle  x_{1}x_{16}-x_{4}x_{17}+\frac{1}{3}x_{6}x_{25}-\frac{2}{3}x_{6}x_{26}-x_{11}x_{14}+x_{9}x_{15} \\
 &&\textstyle  x_{1}x_{18}-x_{4}x_{19}+\frac{1}{3}x_{8}x_{25}-\frac{2}{3}x_{8}x_{26}+x_{11}x_{13}-x_{7}x_{15} \\
 &&\textstyle  x_{1}x_{20}-x_{4}x_{21}-\frac{1}{3}x_{10}x_{25}+\frac{2}{3}x_{10}x_{26}-x_{9}x_{13}+x_{7}x_{14}  \\
&&\textstyle  x_{1}x_{22}-x_{6}x_{19}+x_{8}x_{17}-\frac{1}{3}x_{11}x_{25}-\frac{1}{3}x_{11}x_{26}+x_{5}x_{15} \\
 &&\textstyle  x_{1}x_{23}-x_{6}x_{21}-x_{10}x_{17}+\frac{1}{3}x_{9}x_{25}+\frac{1}{3}x_{9}x_{26}-x_{5}x_{14}  \\
&&\textstyle  x_{1}x_{24}-x_{8}x_{21}-x_{10}x_{19}-\frac{1}{3}x_{7}x_{25}-\frac{1}{3}x_{7}x_{26}+x_{5}x_{13} \\
&&\textstyle  x_{1}x_{2}-x_{11}x_{21}-x_{9}x_{19}-x_{7}x_{17}-\frac{1}{3}x_{5}x_{25}+\frac{2}{3}x_{5}x_{26}\\  
&&\textstyle  x_{1}x_{3}+x_{6}x_{24}-x_{8}x_{23}-x_{10}x_{22}-x_{15}x_{21}-x_{14}x_{19}+\notag\\&&x_{5}x_{12}-x_{13}x_{17}-\frac{1}{3}x_{25}^2+\frac{1}{3}x_{26}^2\\  
&&\textstyle  x_{4}x_{22}-x_{6}x_{18}+x_{8}x_{16}-x_{11}x_{12}+\frac{2}{3}x_{15}x_{25}-\frac{1}{3}x_{15}x_{26}  \\
&&\textstyle  x_{4}x_{23}-x_{6}x_{20}-x_{10}x_{16}+x_{9}x_{12}-\frac{2}{3}x_{14}x_{25}+\frac{1}{3}x_{14}x_{26} \\
 &&\textstyle  x_{4}x_{24}-x_{8}x_{20}-x_{10}x_{18}-x_{7}x_{12}+\frac{2}{3}x_{13}x_{25}-\frac{1}{3}x_{13}x_{26}  \\
&&\textstyle  x_{2}x_{4}-x_{6}x_{24}+x_{8}x_{23}+x_{10}x_{22}-x_{11}x_{20}-x_{9}x_{18}-\notag\\&&x_{7}x_{16}-x_{5}x_{12}+\frac{2}{3}x_{25}x_{26}-\frac{1}{3}x_{26}^2\\  
&&\textstyle  x_{3}x_{4}-x_{15}x_{20}-x_{14}x_{18}-x_{13}x_{16}-\frac{1}{3}x_{12}x_{25}+\frac{2}{3}x_{12}x_{26} \\
 &&\textstyle  x_{6}x_{2}-x_{11}x_{23}-x_{9}x_{22}+x_{5}x_{16}-\frac{2}{3}x_{17}x_{25}+\frac{1}{3}x_{17}x_{26}  \\
&&\textstyle  x_{3}x_{6}-x_{15}x_{23}-x_{14}x_{22}+\frac{1}{3}x_{16}x_{25}+\frac{1}{3}x_{16}x_{26}-x_{12}x_{17}  \\
&&\textstyle  x_{2}x_{8}-x_{11}x_{24}+x_{7}x_{22}+x_{5}x_{18}-\frac{2}{3}x_{19}x_{25}+\frac{1}{3}x_{19}x_{26} \\
 &&\textstyle  x_{3}x_{8}-x_{15}x_{24}+x_{13}x_{22}+\frac{1}{3}x_{18}x_{25}+\frac{1}{3}x_{18}x_{26}-x_{12}x_{19}\\  
&&\textstyle  x_{2}x_{10}-x_{9}x_{24}-x_{7}x_{23}-x_{5}x_{20}+\frac{2}{3}x_{21}x_{25}-\frac{1}{3}x_{21}x_{26} \\
 &&\textstyle  x_{3}x_{10}-x_{14}x_{24}-x_{13}x_{23}-\frac{1}{3}x_{20}x_{25}-\frac{1}{3}x_{20}x_{26}+x_{12}x_{21}\\ 
 &&\textstyle  x_{3}x_{11}-x_{2}x_{15}-\frac{1}{3}x_{22}x_{25}+\frac{2}{3}x_{22}x_{26}+x_{17}x_{18}-x_{16}x_{19}\\ 
 &&\textstyle  x_{3}x_{9}-x_{2}x_{14}+\frac{1}{3}x_{23}x_{25}-\frac{2}{3}x_{23}x_{26}-x_{17}x_{20}+x_{16}x_{21} \\
&&\textstyle  x_{3}x_{7}-x_{2}x_{13}-\frac{1}{3}x_{24}x_{25}+\frac{2}{3}x_{24}x_{26}+x_{19}x_{20}-x_{18}x_{21}  \\
&&\textstyle  x_{3}x_{5}-\frac{1}{3}x_{2}x_{25}-\frac{1}{3}x_{2}x_{26}+x_{17}x_{24}-x_{19}x_{23}+x_{21}x_{22}  \end{eqnarray} \begin{eqnarray} 
&&\textstyle  x_{3}x_{25}-\frac{1}{2}x_{3}x_{26}-\frac{3}{2}x_{2}x_{12}+\frac{3}{2}x_{16}x_{24}-\frac{3}{2}x_{18}x_{23}+\frac{3}{2}x_{20}x_{22}\\ 
&&\textstyle  x_{1}x_{3}-x_{2}x_{4}-x_{6}x_{24}+x_{8}x_{23}+x_{10}x_{22}+x_{11}x_{20}+x_{9}x_{18}-\notag\\&&x_{15}x_{21}+x_{7}x_{16}-x_{14}x_{19}-x_{5}x_{12}-x_{13}x_{17}+\frac{1}{3}x_{25}^2-\frac{1}{3}x_{25}x_{26}+\frac{1}{3}x_{26}^2  
\end{eqnarray}
 \section{SAGE code for the Hilbert series of \(w\SiF\)}\label{sagecode}
 \begin {verbatim}
F4=WeylCharacterRing(['F',4])
W= WeylGroup(['F',4])
L=W.domain()                            
lam=vector([1,0,0,0])                   
a,b,c,d,t,u= var('a,b,c,d,t,u')                         
rho=vector([11/2,5/2,3/2,1/2])
mu=vector([a,b,c,d])
N=[(-1)**(w.length())*(t**(((w*rho)*mu))/ (1-t**(((w*lam)*mu)+u))) for w in W]
num=sum(N)
D=[(-1)**(w.length())*t**(((w*rho)*mu)) for w in W]
den=sum(D)
HilbertSeries=num/den
\end{verbatim}
\bibliographystyle{amsalpha}
\bibliography{imran}

\providecommand{\bysame}{\leavevmode\hbox to3em{\hrulefill}\thinspace}
\providecommand{\MR}{\relax\ifhmode\unskip\space\fi MR }
\providecommand{\MRhref}[2]{%
  \href{http://www.ams.org/mathscinet-getitem?mr=#1}{#2}
}
\providecommand{\href}[2]{#2}
\begin{thebibliography}{MAAvLL92}

\bibitem[AM69]{atiyah}
M.~F. Atiyah and I.~G. Macdonald, \emph{Introduction to commutative algebra},
  Addison-Wesley Publishing Co., Reading, Mass.-London-Don Mills, Ont., 1969.

\bibitem[BE77]{BE}
D.~A. Buchsbaum and D.~Eisenbud, \emph{Algebra structures for finite free
  resolutions, and some structure theorems for ideals of codimension 3}, Amer.
  J. Math \textbf{99} (1977), 447--485.

\bibitem[BKQ]{BKQ}
G.D. Brown, A.~M. Kasprzyk, and M.~I. Qureshi, \emph{Fano manifold in
  {G}orenstein formats}, In preparation.

\bibitem[BKR12]{tom-jerry-1}
G.D. Brown, M.~Kerber, and M.A. Reid, \emph{{Fano 3-folds in codimension 4, Tom
  and Jerry. I.}}, Compos. Math. \textbf{148} (2012), no.~4, 1171--1194.

\bibitem[BKZ]{formats}
G.~Brown, A.~M. Kasprzyk, and L.~Zhu, \emph{Gorenstein formats, canonical and
  {C}alabi-{Y}au threefolds}, arXiv:1409.4644.

\bibitem[Bou02]{bourbaki46}
Nicolas Bourbaki, \emph{Elements of {M}athematics, {L}ie {G}roups and {L}ie
  {A}lgebras, {C}hapters 4--6}, Springer-Verlag, 2002.

\bibitem[CR02]{wg}
A.~Corti and M.~Reid, \emph{Weighted {G}rassmannians}, Algebraic geometry
  (M.~C. Beltrametti, F.~Catanese, C.~Ciliberto, A.~Lanteri, and C.~Pedrini,
  eds.), de Gruyter, Berlin, 2002, pp.~141--163.

\bibitem[dG01]{degraaf}
W.~A. de~Graaf, \emph{Constructing representations of split semisimple {L}ie
  algebras}, J. Pure Appl. Algebra \textbf{164} (2001), no.~1-2, 87--107,
  Effective methods in algebraic geometry (Bath, 2000).

\bibitem[FH91]{harris}
W.~Fulton and J.~Harris, \emph{Representation theory, a first course}, Graduate
  Text in Mathematics, 129, Springer-Verlag, 1991.

\bibitem[GKR07]{rudakov}
A.~L. Gorodentsev, A.~S. Khoroshkin, and A.~N. Rudakov, \emph{On syzygies of
  highest weight orbits}, Moscow {S}eminar on {M}athematical {P}hysics. {II}
  (V.~I. Arnold, D.G. Gindikin, and V.~P. Maslov, eds.), Amer. Math. Soc.
  Transl. Ser. 2, vol. 221, AMS, Providence, RI, 2007, pp.~79--120.

\bibitem[MAAvLL92]{LiE}
A.~M.~Cohen M.~A. A.~van Leeuwen and B.~Lisser, \emph{"{L}i{E}, a {P}ackage for
  {L}ie {G}roup {C}omputations"}, 1992, Computer Algebra Nederland, Amsterdam,
  ISBN 90-74116-02-7.

\bibitem[QS11]{qs}
M.~I. Qureshi and B.~Szendr{\H o}i, \emph{Constructing projective varieties in
  weighted flag varieties}, Bull. Lon. Math Soc. \textbf{43} (2011), no.~2,
  786--798.

\bibitem[QS12]{qs-ahep}
M.~I. Qureshi and Bal{\'a}zs Szendr{\H{o}}i, \emph{Calabi-{Y}au threefolds in
  weighted flag varieties}, Adv. High Energy Phys. (2012), Art. ID 547317, 14
  pp.

\bibitem[Qur15]{qs2}
M.~I. Qureshi, \emph{Constructing projective varieties in weighted flag
  varieties { II}}, Math. Proc. Camb. Phil. Soc. \textbf{158} (2015), 193--209.

\bibitem[Qur17]{QJSC}
Muhammad~Imran Qureshi, \emph{Computing isolated orbifolds in weighted flag
  varieties}, Journal of Symbolic Computation \textbf{79, Part 2} (2017), 457
  -- 474.

\bibitem[Rei87]{YPG}
M.~Reid, \emph{Young person's guide to canonical singularities}, Algebraic
  geometry, {B}owdoin, 1985 ({B}runswick, {M}aine, 1985), Proc. Sympos. Pure
  Math., vol.~46, Amer. Math. Soc., Providence, RI, 1987, pp.~345--414.

\bibitem[Rei15]{codim4}
\bysame, \emph{Gorenstein in codimension 4 - the general structure theory},
  Algebraic Geometry in East Asia (Nov 2011), Advanced Studies in Pure
  Mathematics, vol.~65, Taipei, 2015, pp.~201--227.

\bibitem[Sag]{sage}
\emph{{S}agemath, the {S}age {M}athematics {S}oftware {S}ystem ({V}ersion
  7.5.1)}, {\tt http://www.sagemath.org}.

\bibitem[Ser]{serre}
J.~P. Serre, \emph{Sur les modules projectifs}, Seminaire Dubreil, 1960/61.

\end{thebibliography}

\end{document}